\input amstex
\documentstyle{amsppt}
\magnification 1200
\NoBlackBoxes 
\NoRunningHeads 
\TagsOnRight

\document 
\def\pr{\operatorname{pr}}
\def\DSch{\operatorname{DSch}}

\def\DAlg{\operatorname{DAlg}}
\def\Spec{\operatorname{Spec}}

\def\sep{\operatorname{sep}}
\def\Gal{\operatorname{Gal}}
\def\Hom{\operatorname{Hom}}

\def\id{\operatorname{id}}

\def\Ker{\operatorname{Ker}}
\def\End{\operatorname{End}}
\def\Ann{\operatorname{Ann}}

\def\Ext{\operatorname{Ext}}
\def\FGr{\operatorname{FGr}}
\def\Gr{\operatorname{Gr}}
\def\Mod{\operatorname{Mod}}
\def\Sym{\operatorname{Sym}}

\def\rk{\operatorname{rk}}
\def\det{\operatorname{det}}
\def\Im{\operatorname{Im}}
\def\et{\operatorname{et}}
\def\loc{\operatorname{loc}}

\def\Map{\operatorname{Map}}
\def\Gr{\operatorname{Gr}}
\def\rk{\operatorname{rk}}

\def\det{\operatorname{det}}
\def\GL {\operatorname{GL}}

\def\Frac{\operatorname{Frac}}
\def\SH{\operatorname{SH}}
\def\BR{\operatorname{BR}}

\topmatter
\title Galois Modules arising from Faltings's strict modules
\endtitle
\author Victor Abrashkin
\endauthor
\address Dept.of Maths., Univ. of Durham, Sci. Laboratories, Durham, DH1 3LE, UK 
%\newline  
%Steklov Math. Institute, Gubkina 8, 117 966, Moscow ,Russia 
\endaddress
\email victor.abrashkin{\@ }dur.ac.uk 
\endemail
\date 
\enddate 
\abstract 
Suppose $O$ is a complete discrete valuation 
ring of positive characteristic with perfect residue field. 
The category of finite flat  strict modules was recently 
introduced by Faltings and appears as an equal 
characteristic analogue of the classical category 
of finite flat group schemes. In this paper we obtain 
a classification of these modules and apply it to 
prove analogues of properties, which were 
known earlier for group schemes. 
\endabstract 
\endtopmatter

\subhead 0. Introduction 
\endsubhead 
\medskip 

Throughout all this paper $O$ is the valuation ring of a complete 
discrete valuation field $K$ with perfect 
residue field $k$ of characteristic $p>0$,  
$\Gamma _K$ --- the absolute Galois group 
of $K$ and $\pi $ --- a uniformising element of $K$. 

If characteristic of $K$ is 0, denote by 
$\FGr '(\Bbb Z_p) _O$ the category of finite flat commutative 
group schemes $G$ over $O$ such that the order $|G|$ is a power of $p$. 
Any such group scheme appears as a kernel 
of an isogeny of abelian schemes defined over $O$ 
and reflects important properties of these abelian schemes. 

The classification of objects of the category 
$\FGr '(\Bbb Z_p)_O$ was  
done in [Fo1] under the restriction $e=1$ on the 
absolute ramification index  
$e=e(K)$ of the extension $K/\Bbb Q_p$ 
in terms of finite Honda systems (this classification 
was not complete for $p=2$, for an improved version cf. [Ab1]). 
Further progress was done in papers [Ab3] for $e\leqslant p-1$ 
(group schemes killed by $p$), 
[Co] for $e<p-1$ and, finally, in [Br] for an arbitrary $e$. 

Most interesting number theoretic application of the 
theory of finite flat group schemes 
come from the study of the structure of the 
$\Gamma _K$-module 
$H=G(K_{\sep })$ of geometric points of $G\in\FGr '(\Bbb Z_p)_O$. 
We mention the following three results:
\medskip 

A. {\it Serre's Conjecture} (proved in [Ra]).  
\medskip 

This result describes the action of the inertia subgroup $I_K\subset\Gamma _K$ 
on the semi-simple envelope of $H$. It is given by characters 
$\chi :I_K\longrightarrow k^*$ such that 
for some $N\in\Bbb N$, 
$\chi =\chi _N^a$, where $\chi _N(\tau )=\tau (\pi _N)/\pi _N$, 
$\pi _N^{p^N-1}=\pi $, and 
$a=a_0+a_1p+\dots +a_{N-1}p^{N-1}$ with $p$-digits 
$a_0,a_1,\dots ,a_{N-1}\in [0,e]$;
\medskip 

B.  {\it Ramification estimates}. 
\medskip 

If $p^M\id _G=0$ then the ramification subgroups 
$\Gamma _K^{(v)}$ of $\Gamma _K$ act trivially on $H$ if 
$v>e(M-1+1/(p-1))$, cf. [Fo2];
\medskip 

C. {\it Complete description of $\Gamma _K$-module $H$} (in the case $e=1$, $p>3$ and $pH=0$);
\medskip 

If $e=1$, $p>3$ and $\Bbb F_p[\Gamma _K]$-module $H$ 
satisfies the above Serre's Conjecture and the ramification estimates 
(i.e. the ramification subgroups $\Gamma _K^{(v)}$ act trivially on $H$ 
if $v>\dsize\frac{1}{p-1}$), then 
there is an $G\in\FGr '(\Bbb Z_p)_O$ such that $H=G(K_{\sep })$, cf.[Ab3]. 
\medskip 

Suppose now that characteristic of $K$ is $p$. In this case a 
reasonable analogue of the concept of 
finite flat group scheme 
would give a way to study kernels of isogenies of Drinfeld modules. 
This analogue should appear as a finite flat commutative group scheme $G$  
with continuous  action of a closed subring $O_0=\Bbb F_q[[\pi _0]]\subset O$. 
The notion of $O_0$-module scheme was not very helpful 
until Faltings introduced 
in [Fa] a concept of strict $O_0$-action. The main idea 
of Faltings's definition can be explained as follows. 

Suppose $G=\Spec A$ is a finite flat $O_0$-module over $O$. 
Present $A$ as a quotient 
$O[X_1,\dots ,X_n]/I$ of a ring of polynomials by an ideal 
$I$ and define a deformation $A^\flat $ as $O[X_1,\dots ,X_n]/(I\cdot I_0)$, 
where $I_0=(X_1,\dots ,X_n)$. Faltings requires that $O_0$-module 
structure on $G$, which is given by endomorphisms 
$[o]:A\longrightarrow A$, $o\in O_0$, should have an extension  
to an $O_0$-module structure of the deformation 
$(A,A^\flat )$ of the algebra $A$  
and this extension must satisfy the condition 
of strictness. This means that if 
$[o]^\flat :A^\flat\longrightarrow A^\flat $ 
is an extension of $[o]$ for $o\in O_0$, 
then $[o]^\flat $ must induce multiplication by $o$ 
on $I_0/I_0^2$ and $I/(I\cdot I_0)$. 
This definition gives the category 
$\FGr '(O_0)_O$ of strict $O_0$-modules and its 
objects have many interesting properties discussed in [Fa]. 

In this paper we study number theoretic properties of 
strict $O_0$-modules. In n.1 we present a concept of strict $O_0$-module 
in a slightly different but an equivalent to the original definition by Faltings way.  
In n.2 we describe the category of 
strict $\Bbb F_q$-modules over $O$ and apply this in n.3 
to the classification of objects of the category $\FGr '(O_0)_O$. 
This classification does not depend on the ramification index 
of $K$ over $\Frac O_0=K_0$ and requires only the study of 
primitive elements 
(i.e. the elements $a\in A(G)$ such that $\Delta a=a\otimes 1+1\otimes a$ where 
$\Delta $ is comultiplication) 
of the $O$-algebra $A(G)$ of $G\in\FGr '(O_0)_O$. 
We apply then this classification 
to prove that any object of $\FGr '(O_0)_O$ can be embedded 
into a $\pi _0$-divisible group over $O$. This section contains also 
a comparison of our antiequivalence with parallel results 
in the theory of finite flat group schemes and $p$-adic representations in 
the mixed characteristic case from papers 
[Ab1,3], [Br] and [Fo4]. In n.4 we establish  
precise analogues of the above properties $A,B$ and $C$ 
in the category $\FGr '(O_0)_O$.

The author expresses his gratitude for hospitality to the Max-Planck-Institute 
in Mathematics where a part of this paper was being written. 
\medskip

\subhead 1. Definition and simplest properties 
\endsubhead 
\medskip 

Let $O$ be the valuation ring of a complete discrete valuation field $K$ 
with perfect residue field $k$ of characteristic $p>0$.
All $O$-algebras are assumed usually to be finite, 
i.e. to be free $O$-modules of finite rank. 
\medskip 

\subsubhead {\rm 1.1.} Deformations of augmented $O$-algebras 
\endsubsubhead 
\medskip 

For an augmented $O$-algebra $A$, we agree to use the following notation: 
$\varepsilon _A:A\longrightarrow O$ --- the morphism of augmentation;  
$\Ker\varepsilon _A=I_A$ --- the augmentation ideal; 
$A^{loc}=\mathbin{\underset n \to \varprojlim}A/I_A^n$ 
--- the $I_A$-completion of $A$; $\eta _A:A\longrightarrow A^{loc}$ --- a natural projection. 
Notice,  the correspondence $A\mapsto A^{loc}$ is functorial.   

The objects of the category $\DAlg _O$ are the triples $\Cal A=(A,A^\flat ,i_{\Cal A})$, where 
$A$ is a finite $O$-algebra, $A^\flat $ is an $O$-algebra such that there is  
a polynomial ring  
$O[\bar X]=O[X_1,\dots ,X_n]$, $n\geqslant 0$, and its ideal $I$ such that    
$A^{loc}=O[\bar X]/I$ and 
\linebreak 
$A^\flat =O[\bar X]/(I\centerdot I_0)$ with 
$I_0:=(X_1,\dots ,X_n)\supset I$, and 
$i_{\Cal A}:A^\flat \longrightarrow A^{loc}$ is a 
natural epimorphism of $O$-algebras. 

A morphism $\bar f=(f,f^\flat ):\Cal A\longrightarrow \Cal B=(B,B^\flat ,i_{\Cal B})$ 
in $\DAlg _O$ is given by a morphism of augmented  
$O$-algebras $f:A\rightarrow B$ and an $O$-algebra 
morphism $f^\flat :A^\flat \longrightarrow B^\flat $ 
such that $i_{\Cal A}\circ f^{loc}=f^\flat \circ i_{\Cal B}$.

In the category $\DAlg _O$, 
$\Cal O=(O,O,\id _O)$ is an initial object and any 
$\Cal A=(A,A^{\flat }, i_{\Cal A})$ has a natural augmentation to $\Cal O$, 
$\varepsilon _{\Cal A}=(\varepsilon _A, \varepsilon _{A^\flat }):\Cal A\longrightarrow\Cal O$, 
where $\Ker\varepsilon _{A^\flat }=\Ann (\Ker i_{\Cal A}):=I_{A^\flat }$.

Notice that the $O$-modules  
$t^* _{\Cal A}=I_0/I_0^2$ (this module is free) and $N_{\Cal A}=I/(I\cdot I_0)$ 
do not depend on the choice of the covering $O[\bar X]\longrightarrow A^{\flat }$.  
Indeed, the first coincides with 
$I_{A^\flat }/I^2_{A^\flat }$ and the second --- with $\Ker i_{\Cal A}$.

If $\Cal A=(A,A^\flat ,i_{\Cal A})$ and $\Cal B=(B,B^\flat ,i_{\Cal B})$ 
are objects of $\DAlg _O$ and $f:A\longrightarrow B$ is a morphism of augmented 
$O$-algebras, then the set of all $f^\flat $ such that 
$(f,f^\flat )\in\Hom _{\DAlg _O}(\Cal A,\Cal B)$ is not empty and has  
a natural structure of a principal homogeneous space over the group 
$\operatorname{Hom}_{O\operatorname{-mod}}(t^* _{\Cal A},N_{\Cal B})$.  

The deformation 
$\Cal A=(A,A^\flat ,i_{\Cal A})$ of $A$ will be called minimal if $i_{\Cal A}$ 
induces isomorphism of $t^* _{\Cal A}\otimes k$ onto 
$I_A/I^2_A\otimes k$. In other words, $f$ is minimal if minimal 
systems of generators 
of $O$-algebras $A^\flat $ and $A^{\loc }$ contain 
the same number of elements. It is easy to see that if $\Cal A$ is minimal, 
$\Cal A'=(A,{A'}^\flat ,i_{\Cal A'})\in\DAlg _O$ and 
$\bar f=(\id _A,f^\flat )\in\Hom _{\DAlg _O}(\Cal A,\Cal A')$ 
then there is an $\bar g\in\Hom _{\DAlg _O}(\Cal A',\Cal A)$  such that 
$\bar f\circ\bar g=\id _{\Cal A}$. 
In particular, all minimal deformations of a given $O$-algebra $A$ are isomorphic in $\DAlg _O$. 

Let $\DAlg ^*_R$ be a quotient category for $\DAlg _O$: it has the same 
objects but its morphisms 
are equivalence classes of morphisms from $\Hom _{\DAlg _O}(\Cal A,\Cal B)$  
arising  from the same $O$-algebra morphisms $f:A\longrightarrow B$. Then the forgetful functor 
$\Cal A=(A,A^\flat ,i_{\Cal A})\mapsto A$ is an equivalence of $\DAlg ^*_O$ and the category of 
augmented finite $O$-algebras. 
\medskip 

\subsubhead {\rm 1.2.} Deformations of affine group schemes 
\endsubsubhead
\medskip 

Let $\DSch _O$ be the dual category for $\DAlg _O$. Its objects 
appear in the form $\Cal H=\Spec\Cal A=(H,H^\flat ,i_{\Cal H})$, where 
$H=\Spec A$ and $H^\flat =\Spec A^\flat $ are finite flat $O$-schemes, 
$\Cal A=(A,A^\flat ,i_{\Cal A})\in\DAlg _O$,  
and $i_{\Cal H}:H\rightarrow H^\flat $ 
is a closed embedding of $O$-schemes. 
This category has direct products:  
if for $i=1,2$, $\Cal A_i=(A_i,A_i^\flat ,i_{\Cal A_i})$ with 
$A^{loc}_i=O[\bar X_i]/I$, $A_i^\flat =O[\bar X_i]/(I_i\cdot I_{0i})$, then 
the product $\Spec\Cal A_1\times\Spec\Cal A_2$ is given by $\Spec (\Cal A_1\otimes\Cal A_2)$, where 
 $\Cal A_1\otimes\Cal A_2:=(A_1\otimes _OA_2,(A_1\otimes _OA_2)^\flat ,\kappa )$, 
$(A_1\otimes _OA_2)^\flat $ is the quotient of $O[\bar X_1\otimes 1, 1\otimes \bar X_2]$ by 
the product of ideals $I_1\otimes 1+1\otimes I_2$ and $I_{01}\otimes 1+1\otimes I_{02}$ and 
$\kappa $ is the natural projection. Notice that for $i=1,2$, 
two projections $\pr _i$ from this product to its components $\Spec\Cal A_i$ come from the 
natural embeddings of $O[\bar X_i]$ into $O[\bar X_1\otimes 1,1\otimes\bar X_2]$.

Let $\FGr _O$ be the category of group objects in $\DSch _R$. If $\Cal G=\Spec\Cal A\in\FGr _O$  
then its group structure is given via comultiplication  
$\bar\Delta =(\Delta ,\Delta ^\flat ):\Cal A\longrightarrow\Cal A\otimes\Cal A$, 
counit $\bar\varepsilon =(\varepsilon, \varepsilon ^\flat ):\Cal A\longrightarrow\Cal O$ and  
coinversion $\bar\imath =(i,i^\flat ):\Cal A\longrightarrow\Cal A$ morphisms, which  satisfy 
usual axioms. 
The morphisms in $\FGr _O$ are morphisms of group objects. 
As usually, $\FGr _O$ is an additive category.

Notice that 
\newline 
a) $G=\Spec A$ is a finite flat group scheme over $O$ with comultipilication $\Delta $, 
counit $\varepsilon $ and coinversion $i$;
\newline 
b) $\bar\varepsilon =\varepsilon _{\Cal A}$, cf. n.1.1. 
\newline 
c) the counit axiom gives for $i=1,2$, $\Delta _i^\flat\circ\pr _i=\id _{A^\flat }$ 
and implies a uniqueness of $\Delta ^\flat$ as a lifting of $\Delta $; 
\newline 
d) if $\Cal A=(A,A^{\flat}, i_{\Cal A})\in\DAlg _O$ and $G=\Spec A$ is an affine group scheme then 
there is a unique structure of group object on $\Spec\Cal A$ compatible with that of $G$;   
\newline 
e) if 
$f:G\longrightarrow H$ is a morphism of group schemes and 
$(f,f^\flat )\in\Hom _{\DSch _O}(\Cal G, \Cal H)$ then 
$(f,f^\flat )\in\Hom _{\FGr _O}(\Cal G,\Cal H)$. 
\medskip 

The above properties have the following 
interpretation. 
Define the quotient category 
$\FGr ^*_O$ as the category consisiting of the objects of the category 
$\FGr _O$ but where $\Hom _{\FGr ^*_O}(\Cal G,\Cal H)$ consists of equivalence classes 
of morphisms from the category $\FGr _O$ which induce the same morphisms 
of affine group schemes $G\rightarrow H$. Then the natural functor 
$\Cal G\mapsto G$ is an equivalence of categories. 
We can use this equivalence to define the abelian groups of 
equivalence classes of 
short exact sequences $\Ext (\Cal G, \Cal H)$ in $\FGr _O$. 
These groups are functorial by both arguments and there are standard 
6-terms exact $\Hom -\Ext $-sequences. 
\medskip

\subsubhead{\rm 1.3.} The categories of strict $R$-modules 
\endsubsubhead 
\medskip

Suppose $R$ is a ring and  $O$ is an $R$-algebra. 

Suppose $\Cal G$ is an $R$-module object in the category $\DSch _O$.
Then $\Cal G$ is an object of $\FGr _O$ and there is a map 
$R\longrightarrow\End _{\FGr _O}(\Cal G)$ satisfying the  usual axioms from the definition of 
$R$-module. For $r\in R$ and $\Cal G=\Spec\Cal A$, denote by 
$\bar {[r]}=([r],[r]^\flat )$ the morphism of action of $r$ on 
$\Cal A=(A,A^\flat ,i_{\Cal A})$. Clearly, $G=\Spec A$ is an $R$-module in the 
category of finite flat schemes over $O$. For any such $G$, 
the $R$-module structure on the deformation 
$(G,G^\flat ,i)\in\FGr _O$ is given by liftings 
$[r]^\flat :A^\flat\longrightarrow A^\flat $ of morphisms $[r]:A\longrightarrow A$, $r\in R$. 
Notice that $[r]^\flat $ are morphisms of augmented algebras. 
All such liftings are automatically compatible with the group structure on 
this deformation, i.e. for any $r\in R$, 
it holds $[r]\circ\Delta ^\flat =\Delta \circ ([r]\otimes [r])$. 
So, the above system gives an $R$-module 
structure if and only if for any $r_1,r_2\in R$, 
$$[r_1+r_2]^\flat =\Delta ^\flat\circ ([r_1]\otimes [r_2])^\flat ,\ \ \ 
[r_1r_2]^\flat =[r_1]\circ [r_2]^\flat \tag{1}$$
where $([r_1]\otimes [r_2])^\flat $ is induced by $[r_1]^\flat\otimes [r_2]^\flat $. 

We denote by $\FGr (R)_O$ the category of such $R$-module objects 
$\Cal G$ where $R$-action is strict, i.e. if $\Cal G=\Spec\Cal A$ 
then any $r\in R$ acts on $t^* _{\Cal A}$ and $N_{\Cal A}$ via 
scalar multiplication by $r$. This is a basic definition from Faltings' paper [Fa]. 

Suppose $\Cal G_1=(G,G_1^\flat ,i_1)\in\FGr _O$ and 
$\Cal G_2=(G,G_2^\flat ,i_2)\in\FGr _O$ 
are two deformations of a finite flat group scheme $G$ over $O$. 
By 1.2, $\Cal G_1$ and $\Cal G_2$ are isomorphic 
in the category $\FGr _O^*$. Suppose $\Cal G_1$ is equiped with the strict $R$-action. 
Then there is a unique (strict) $R$-action 
on $\Cal G_2$ such that any $(\id _G,\phi )\in\Hom_{\FGr _O}(\Cal G_1,\Cal G_2)$ 
and any $(\id _G,\psi )\in\Hom_{\FGr _O}(\Cal G_2,\Cal G_1)$ are, actually, 
morphisms in the category $\FGr (R)_O$.

Denote by $\FGr^*(R)_O$ the quotient category of $\FGr (R)_O$ where morphisms  
 are equivalence classes 
of morphisms 
$(G,G^\flat ,i)\longrightarrow (H,H^\flat ,j)$ 
in the category $\FGr (R)_O$ inducing the same morphism $G\longrightarrow H$. 
By the above property, all isomorphism classes of objects in $\FGr (R)_O$ appear as  
$R$-module finite flat schemes $G$ together with a lifting of its 
$R$-action to some chosen deformation $G^\flat $ which satisfies the above conditions (1).

For example, if $q=p^n$ with $n\in\Bbb N$, the objects of the category 
$\FGr (\Bbb F_q)_O$ appear as $\Spec\Cal A$, where 
$\Cal A=(A,A^\flat ,i_{\Cal A})$, 
$A^{loc}=O[\bar X]/I$ and $A^\flat =O[\bar X]/(I\cdot I_0)$, and $\Bbb F_q$-action 
is induced by 
the scalar action of $\Bbb F_q$ on $O[\bar X]$ 
(i.e. $[\alpha ](\bar X)=\alpha \bar X$, $\alpha\in\Bbb F_q$) 
and there are generators $j_1,\dots ,j_n$ of the ideal $I$ such that 
$[\alpha ]j_i=\alpha j_i$ for all $i=1,\dots ,n$ and $\alpha\in\Bbb F_q$. 
 
If $R=\Bbb F_q[\pi ]$ with an indeterminate $\pi $ and $\Cal G\in\FGr (R)_O$ then 
$\Cal G\in\FGr (\Bbb F_q)_O$ and (in addition to the above assumptions) 
$R$-action will be determined completely by the action of $\pi $ given by 
the correspondence 
$$\bar X\mapsto [\pi ](\bar X)=\pi\bar X+\bar F(\bar X)$$
where $\bar F$ is any vector power series with coefficients in $O$ from $I_0^2$ 
such 
that $\bar F(\alpha\bar X)=\alpha\bar F(\bar X)$ for all $\alpha\in\Bbb F_q$. 
This action is strict iff $[\pi ]j_i\equiv\pi j_i\operatorname{mod}(I\cdot I_0)$ 
for the above generators $j_1,\dots ,j_n$ of $I$. 
\medskip 

1.4. Suppose $i:\Cal G_1\longrightarrow \Cal G$ is a 
closed embedding of strict $R$-modules. Then the quotient 
$\Cal G_2=\Cal G/\Cal G_1$ has a natural structure of a 
strict $R$-module scheme and the projection 
$j:\Cal G\longrightarrow\Cal G_2$ is a morphism in 
the category $\FGr (R)_O$. 
Similarly, if $j:\Cal G\longrightarrow\Cal G_2$ is 
a fully faithful morphism of strict $R$-module schemes 
then its kernel $\Cal G_1$ is a strict $R$-module and 
its (closed) embedding $i:\Cal G_1\longrightarrow\Cal G$ 
is a morphism from $\FGr (R)_O$.

Notice that for any $\Cal G=(G,G^\flat ,i_{\Cal G})\in\FGr (R)_O$, 
the short exact sequence of group schemes 
$$0\longrightarrow G^{\loc }\longrightarrow G
\longrightarrow G^{\et }\longrightarrow 0\tag{2}$$
(here,  $G^{\loc }$ is the maximal local subgroup scheme of $G$ and $G^{\et }$ 
is the maximal etale quotient of $G$)  
induces a short exact sequence in $\FGr (R)_O$ 
$$0\longrightarrow\Cal G^{\loc }\longrightarrow
\Cal G\longrightarrow\Cal G^{\et }\longrightarrow 0\tag{3}$$
where $\Cal G^{\loc }=(G^{\loc },G^\flat )$, $\Cal G^{\et }=(G^{\et }, \Spec O)$.

Finally, notice that over the valution ring $O_1$ of some 
unramified extension $K_1$ of $K$, $G^{\et }\otimes O_1$ 
is constant, i.e. 
$A(G^{\et })\otimes O_1=\Map (H^{\et }, O_1)$ with $H^{\et }=G^{\et }(K_1)$. 
This also implies the existence of a finite extension 
$K_1'$ of $K$ such that the short 
exact sequence (2) (as well as (3)) splits over 
the valuation ring of $K_1'$.  
\medskip 

\ \

\subhead 2. Group schemes with strict $\Bbb F_q$-action  
\endsubhead 
\medskip 

\ \ 

As earlier, $O$ is the valuation ring of a complete discrete valuation field $K$ 
of characteristic $p$  
with perfect residue field $k$, $q=p^{N_0}$ with $N_0\in\Bbb N$. 
In this section we assume that $O$ is an $\Bbb F_q$-algebra 
(then $\Bbb F_q\subset k$ and $O\simeq k[[\pi ]]$, where 
$\pi $ is a uniformiser in $K$)   
and study the full subcategory $\FGr '(\Bbb F_q)_O$ of 
the category $\FGr (\Bbb F_q)_O$ consisting 
of strict finite $\Bbb F_q$-modules $G$ over $O$ 
with etale generic fibre. As it was noticed in n.1, its objects 
appear as finite flat group schemes 
$G=\Spec A(G)$ over $O$ with $\Bbb F_q$-action 
such that $A(G^{loc})$ (or even $A(G)$) can be presented in the form $O[\bar X]/I$ in 
such a way that $\Bbb F_q$-action on $G$ is induced by a scalar $\Bbb F_q$-action on coordinates of 
$\bar X$ and on some system of generators of the ideal $I$.  
\medskip 

\subsubhead {\rm 2.1.} The category $\Mod _{O,\sigma }$ and the functor $\Gr $
\endsubsubhead 
\medskip 

Let $\sigma :O\rightarrow O$ be the endomorphism of $q$-th power. For any $O$-module $M$, 
set $M_{\sigma }=M\otimes _OO$ where an $O$-module structure on $O$ is given via $\sigma $. 
Objects of the category $\Mod _{O,\sigma }$ are free $O$-modules $M$ of finite rank 
with an $O$-linear morphism $\Phi :M_\sigma\rightarrow M$ such that its $K$-linear extension 
$\Phi _K$ is an isomorphism of $K$-vector spaces $M_{\sigma ,K}=
M_\sigma \otimes _OK$ and $M_K=M\otimes _OK$. 
Morphisms in $\Mod _{O,\sigma }$ are morphisms of $O$-modules commuting with 
$\Phi $-action. 

For an $M\in\Mod _{O,\sigma }$, introduce the $O$-algebra 
$$A[M]=\Sym _OM/\{\Phi m-m^{\otimes q}\ |\ m\in M\}$$
 In other words, if 
$m_1,\dots ,m_n$ is an $O$-basis of $M$ and $\Phi m_i=\sum _j o_{ij}m_j$ (where $\det (o_{ij})\ne 0$),  
then $A[M]=O[X_1,\dots ,X_n]/(\psi _1,\dots ,\psi _n)$, where 
$\psi _i=X_i^q-\sum_jo_{ij}X_j$, $i=1,\dots ,n$. This $O$-algebra $A[M]$ 
does not depend on the choice of a basis $m_1,\dots ,m_n$ and gives rise to 
the strict finite flat $\Bbb F_q$-module scheme $G[M]=\Spec A[M]$ with 
comultiplication $\Delta $ such that $\Delta (X_i)=X_i\otimes 1+1\otimes X_i$ and an 
$\Bbb F_q$-action such that $[\alpha ](X_i)=\alpha X_i$, where $\alpha \in\Bbb F_q$, $i=1,\dots ,n$. 
The equality $\Phi (M_{\sigma ,K})=M_K$ implies that the generic fibre of 
$G[M]$ is etale, and therefore, $G[M]\in\FGr '(\Bbb F_q)_O$. 

Identify $M$ with its image in $A[M]$ via the natural map $m_i\mapsto X_i$, $i=1,\dots ,n$.

\proclaim{Proposition 1} 
\newline 
{\rm a)} $M=\{m\in A[M]\ |\ \Delta m=m\otimes 1+1\otimes m,  
\ \ [\alpha ]m=\alpha m, \alpha\in\Bbb F_q\}$; 
\newline 
 {\rm b)} the natural inclusions $M_i\subset A[M_i]$, $i=1,2$, 
induce the identification 
$$\Hom _{\Mod _{O,\sigma }}(M_1,M_2)=\Hom _{\FGr (\Bbb F_q)_O}(G[M_2],G[M_1]).$$
\endproclaim 

This proposition implies that the correspondence $M\mapsto G[M]$ induces a 
fully faithful functor 
$\Gr :\Mod _{O,\sigma }\longrightarrow\FGr '(\Bbb F_q)_O$. 
The rest of this section will be devoted to the proof of 
the following theorem. 

\proclaim{Theorem 1} The functor $\Gr $ is an antiequivalence of categories.
\endproclaim

Notice that the correspondence $\Gr \otimes K: M _K\mapsto G[M](K_{sep})$ is an 
antiequivalence 
of the category of $\sigma $-etale $K$-modules and the category of 
$\Bbb F_q[\Gamma _K]$-modules, cf. [Fo4]. 
\medskip 

2.2. Later we need the following 
property of the Hochschild cohomology of group schemes $G[M]$. 

\proclaim{Lemma 1} With the above notation suppose that 
$a\in A[M]_K$ is such that 
$\delta ^+a:=\Delta(a)-a\otimes 1-1\otimes a\in A[M]\otimes A[M]$ 
and $[\alpha ]a=\alpha a$ for all $\alpha\in\Bbb F_q$. 
Then there is an $a'\in A[M]$ such that $\delta ^+a=\delta ^+a'$ and $a-a'\in M_K$.
\endproclaim 

\demo{Proof of Lemma 1} Let $\delta ^+_k=\delta ^+\otimes k:A[M]_k\longrightarrow 
A[M]_k^{\otimes 2}$, where $A[M]_k=A[M]\otimes k$. It will be sufficient 
to prove that 
$$\{b\in A[M]_k\ |\ \delta ^+_kb=0, \ [\alpha ]b=\alpha b, \ \forall\alpha\in\Bbb F_q\}=M_k\tag{4}$$
Note that 
$$\{X_1^{i_1}\dots X_n^{i_n}\ |\ 0\leqslant i_1,\dots ,i_n<q\}$$
is a $k$-basis of $A[M]_k$ and 
$$\{X_1^{j_1}\dots X_n^{j_n}\otimes X_1^{l_1}\dots X_n^{l_n}\ |
\ 0\leqslant j_1,\dots ,j_n,l_1,\dots ,l_n<q\}$$
is a $k$-basis of $A[M]_k^{\otimes 2}$. Now the equality (4) 
is implied by the following 3 observations:

a) for any $0\leqslant i_1,\dots ,i_n<q$, $\delta ^+_k(X_1^{i_1}\dots X_n^{i_n})$ is 
a linear combination of 
\linebreak 
$X_1^{j_1}\dots X_n^{j_n}\otimes X_1^{l_1}\dots X_n^{l_n}$ with 
$j_1,\dots ,j_n,l_1,\dots ,l_n\geqslant 0$ such that $j_1+l_1=i_1$, 
\linebreak 
$\dots , j_n+l_n=i_n$;  
\medskip 

b) if $\mathbin{\underset{0\leqslant i_1,\dots ,i_n<q}\to\sum }\alpha _{i_1\dots i_n}
X_1^{i_1}\dots X_n^{i_n}\in\Ker\delta ^+_k$ then for each multi index $(i_1,\dots ,i_n)$, 
either 
$\alpha _{i_1\dots i_n}=0$, or  
$\delta ^+_k(X_1^{i_1}\dots X_n^{i_n})=0$; 
\medskip 

c) for any $0\leqslant i_1,\dots ,i_n<q$, 
$\delta ^+_k(X_1^{i_1}\dots X_n^{i_n})=0$ iff 
\newline  
$X_1^{i_1}\dots X_n^{i_n}\in\{X_s^{p^t}\ |\ 1\leqslant s\leqslant n, 0\leqslant t<N_0\}$. 
\enddemo 
\medskip

2.3. For $G\in\FGr (\Bbb F_q)_O$, let 
$$L(G)=\{a\in A(G)\ |\ \Delta a=a\otimes 1+1\otimes a, 
\ [\alpha ]a=\alpha a,  \forall\alpha\in\Bbb F_q\}.$$
Clearly, the above lemma implies that 
if $M\in\Mod _{O,\sigma }$, then 
$L(G[M])=M$. 

The following lemma can be proved directly from the above definition of 
$L(G)$. 

\proclaim{Lemma 2} 
\newline 
{\rm a)} If $O_1$ is a complete discrete valuation ring containing $O$, then 
\newline 
$L(G\otimes O_1)=L(G)\otimes O_1$;
\newline 
{\rm b)} For any two objects $G_1$ and $G_2$ of $\FGr (\Bbb F_q)_O$, it holds 
\newline 
$L(G_1\times G_2)=L(G_1)\oplus L(G_2)$. 
\endproclaim 

The $q$-th power map on $A(G)$ induces 
$\Phi :L(G)_\sigma \longrightarrow L(G)$ which provides 
$L(G)$ with the structure of an object of the category 
$\Mod _{O,\sigma }$. Indeed, 
if $H=G(K_{\sep })$ with the natural structure of 
$\Bbb F_q[\Gamma _K]$-module, then 
$L(G)_K=\Hom _{\Bbb F_q[\Gamma _K]\operatorname{-mod}}(H,K_{\sep })$ and 
$$\Im\Phi _K=\Hom _{\Bbb F_q[\Gamma _K]\operatorname{-mod}}(H,\sigma (K_{\sep }))
\otimes _{\sigma (K)}K=L(G)_K$$ 
because  
$\sigma (K_{\sep })\otimes _{\sigma (K)}K=K_{\sep }$.

 The correspondence $L:G\mapsto L(G)$ 
is the functor from $\FGr '(\Bbb F_q)_O$ to $\Mod _{O,\sigma }$. 
Clearly, the composition $\Gr\circ L$ is equivalent to the identity functor 
on $\Mod _{O,\sigma }$. 
So, $\Gr $ is an antiequivalence (with the quasi-inverse functor $L$) if any 
$G\in\FGr '(\Bbb F_q)_O$ is isomorphic to $G[M]$ for a suitable 
$M\in\Mod _{O,\sigma }$. When proving this property we can enlarge (if necessary) 
the basic ring $O$ by the following lemma.

\proclaim{Lemma 3} If $O'$ is a complete discrete valuation ring containing $O$ 
and $G\in\FGr '(\Bbb F_q)_O$ then $G\simeq G[M]$ with  
$M\in\Mod _{O,\sigma }$ 
if and 
only if $G\otimes O_1\simeq G[M_1]$ for some 
$M_1\in\Mod _{O_1,\sigma }$.
\endproclaim 

\demo{Proof} Clearly, if  $G\simeq G[M]$ with $M\in\Mod _{O,\sigma }$, 
then $G\otimes O_1\simeq G[M\otimes O_1]$. 

If $G\otimes O_1\simeq G[M_1]$ with $M_1\in\Mod _{O,\sigma }$, 
then $M_1\simeq L(G\otimes O_1)=L(G)\otimes O_1$, therefore, 
$G\otimes O_1\simeq G[L(G)]\otimes O_1$ and the natural inclusion 
$G[L(G)]\subset G$ is an isomorphism. Lemma 3 is proved
\enddemo

\subsubhead{\rm 2.4.} $\Bbb F_q$-module schemes of order $q$
\endsubsubhead 

\proclaim{Proposition 2} Suppose $G$ is an $\Bbb F_q$-module scheme of order $q$ 
such that $G\otimes K$ is etale and $\Bbb F_q$ acts on the cotangent space 
$t_G^*$ via scalar multiplication. Then $G=G[M]$, where 
$M=Om$ and $\Phi m=\lambda m$ with some $\lambda\in O$, $\lambda\ne 0$ 
(in particular, $G$ is a strict $\Bbb F_q$-module). 
\endproclaim 

\demo{Proof} Let $I_G=\Ker e_G$ (where $e_G$ is the counit morphism). 
If $I_G=I_G^2$ then $G$ is etale and over some bigger 
(unramified over $O$) basic ring $O_1$, 
$G\otimes O_1$ is constant, i.e. $G\otimes O_1\simeq G[M_1]$, 
where $M_1=O_1m_1$ and $\Phi (m_1)=m_1$. Therefore, $G\simeq G[M]$ where 
$M\otimes O_1\simeq M_1$, i.e. $M=Om$, $m=\xi m_1$ with $\xi\in O_1^*$ and 
$\Phi m=\lambda m$  with $\lambda =\xi ^{q-1}\in O\cap O_1^*=O^*$.

Suppose now that $I_G\ne I_G^2$. Consider the decomposition 
$I_G=\oplus _{1\leqslant n<q}I_{G,n}$, where 
an $\alpha\in\Bbb F_q$ acts on $I_{G,n}$  via multiplication by $\alpha ^n$. 

Because $t_G^*=I_G/I^2_G\ne 0$, $\rk _OI_{G,1}\geqslant 1$ and, therefore, for all 
$1\leqslant n<q$, $\rk _OI_{G,n}\geqslant 1$ (if $x\in I_{G,1}$, $x\ne 0$, 
then $x^n\ne 0$ and $x^n\in I_{G,n}$). So, 
all these ranks are equal to 1 and, if $I_{G,1}=xO$, then the equality 
$I_G=I_{G,1}+I_G^2$ implies for all 
$1\leqslant n<q$, $I_{G,n}=x^nO$. 

Finally, $x^q\in I_{G,1}$ and $x^q=\lambda x$ for some $\lambda\in O$. 
Clearly, $\lambda\ne 0$ because $G\otimes K$ is etale.  
This proves that $A(G)=A[M]$ with $M=Om$, $\Phi m=\lambda m$.  

Suppose $\Delta x=x\otimes 1+1\otimes x+\delta $. Then $\delta ^q=\lambda \delta $ and 
$\delta\in (x\otimes x)A(G)\otimes A(G)$. This implies easily that 
$\delta =0$. It remains to note that by the choice of $x$, 
$[\alpha ]x=\alpha x$ for any $\alpha\in\Bbb F_q$, and this 
action of $\Bbb F_q$ is strict. 
\enddemo 

\definition{Definition} We shall denote by $\mu _{\lambda }$ 
the $\Bbb F_q$-module scheme $G[M]$ from the above proposition.
\enddefinition 

2.5. Suppose $G\in\FGr '(\Bbb F_q)_O$ satisfies the following four assumptions:

1) $G=G^{loc}$;

2) there is a short exact sequence in the category $\FGr '(\Bbb F_q)_O$
$$0\longrightarrow\mu _{\lambda }\mathbin{\overset i\to
\longrightarrow} G\longrightarrow G_1\longrightarrow 0$$
where $\lambda\in O, \lambda\ne 0$;

3) the above exact sequence splits over $K$, i.e. there is a $\Gamma _K$-invariant 
section of the projection of $\Bbb F_q[\Gamma _K]$-modules $G(K_{sep})\longrightarrow G_1(K_{sep})$;

4) there is an $M_1\in\Mod _{O,\sigma }$ such that $G_1\simeq G[M_1]$.
\medskip

\proclaim{Proposition 3} With the above assumptions 1)-4), there is an 
$M\in\Mod _{O,\sigma }$ such that $G\simeq G[M]$.
\endproclaim 

\demo{Proof} Let $A(\mu _{\lambda })=O[u]$ with $u^q=\lambda u$, 
$\Delta u=u\otimes 1+1\otimes u$, $[\alpha ]u=\alpha u$ for $\alpha\in\Bbb F_q$. 
Suppose the closed embedding $i:\mu _{\lambda }\longrightarrow G$ is given by 
the $O$-algebra morphism $i^*:A\longrightarrow A(\mu _{\lambda })$. 
Set $A=A(G)$ and $B=A(G_1)$. Then 
$$s=\Delta _G\circ (i^*\otimes\id _A):A\longrightarrow A(\mu _{\lambda })\otimes A$$
is the coaction of $\mu _{\lambda }$ on $A$ and $B=A^{\mu _{\lambda }}=\{a\in A\ |\ s(a)=1\otimes a\}$. 

Choose an $a_0\in A$ such that $i^*(a_0)=u^{q-1}$ and $[\alpha ]a_0=a_0$ 
for all $\alpha\in\Bbb F_q$. Set 
$$s(a_0)=1\otimes a_0+u\otimes a_1+\dots +u^{q-1}\otimes a_{q-1}$$
where $a_1,\dots ,a_{q-1}\in A$. Then for any $\alpha\in\Bbb F_q$ and 
$m=1,\dots ,q-1$, $[\alpha ]a_m=\alpha ^{-m}a_m$ and the identity 
$$(\id\otimes i^*)(s(a_0))=\Delta (i^*a_0)=\Delta (u^{q-1})=(u\otimes 1+1\otimes u)^{q-1}$$
implies that $i^*(a_{q-1})=1$ and $i^*(a_{q-2})=-u$. Therefore, 
$a_{q-1}\in A^*$ (because $A=A^{loc}$) and $A=B[a_{q-2}]$ (because 
$a_{q-2}$ generates $A$ modulo a nilpotent ideal $I_BA$). 

The equality 
$$1\otimes s(a_0)+u\otimes s(a_1)+\dots +u^{q-1}\otimes s(a_{q-1})=
(\id \otimes s)(s(a_0))$$
$$=(\Delta\otimes\id )(s(a_0))=1\otimes 1\otimes a_0+(u\otimes 1+1\otimes u)\otimes a_1+
\dots +(u\otimes 1+1\otimes u)^{q-1}\otimes a_{q-1}$$
implies that $s(a_{q-1})=1\otimes a_{q-1}$ and $s(a_{q-2})=1\otimes a_{q-2}-u\otimes a_{q-1}$. 

Therefore, $a_{q-1}\in B\cap A^*=B^*$ and for $\theta =-a_{q-2}/a_{q-1}$, we have 
$A=B[\theta ]$, $s(\theta )=1\otimes\theta +u\otimes 1$, $[\alpha ]\theta =\alpha\theta $ for 
$\alpha\in\Bbb F_q$, and $\theta ^q-\lambda\theta =b\in B$. 

Notice also that $\delta ^+\theta =\Delta\theta -\theta\otimes 1-1\otimes\theta \in B\otimes B$, 
because this is a $(\mu_{\lambda }\times\mu_{\lambda })$-invariant element of $A\otimes A$. 

From condition 3) it follows the existence of $v\in A_K$ such that 
$i^*(v)=u$, $v^q=\lambda v$, $\Delta v=v\otimes 1+1\otimes v$ and $[\alpha ]v=\alpha v$ 
for $\alpha\in\Bbb F_q$. Therefore, denoting the  $K$-linear extension of $s$ by $s_K$, we have 
$s_K(v)=1\otimes v+u\otimes 1$. This implies that 
$\theta -v=b_0\in B\otimes K$, 
and $\delta ^+b_0=\delta ^+\theta $. 

Therefore, by the property of Hochshild cohomology from the end of n.2.1,  
there is a $b_1\in B$ (with the property 
$[\alpha ]b_1=\alpha b_1$,  $\forall\alpha\in\Bbb F_q$) 
such that $\delta ^+b_1=\delta ^+\theta $ and replacing $\theta $ by 
$\theta -b_1$ we can assume that $\Delta \theta =\theta\otimes 1+1\otimes\theta $. 
Therefore, $M=M_1+O\theta \in\Mod _{O,\sigma }$, $A=A[M]$ and $G=G[M]$. 
The proposition is proved.  
\enddemo 

2.6. Now we can finish the proof of Theorem 1. 

Suppose $G\in\FGr '(\Bbb F_q)_O$ is of order $q^N$, $N\in\Bbb N$. As it was noticed in n.2.3,  
it will be sufficient to prove that $G\simeq G[M]$ for a suitable $M\in\Mod _{O,\sigma }$. 
When proving this property we can enlarge (if necessary) the basic ring $O$. 

Apply induction on $N$. 

The case $N=1$ follows from Proposition 2 in n.2.4. 

Suppose that $G\in\FGr '(\Bbb F_q)_O$ is of order $q^{N}$ and any 
$\Bbb F_q$-module scheme of order $q^{N-1}$ appears in the form 
$G[M_1]$ with $M_1\in\Mod _{O,\sigma }$. 
By n.1.4 
and Lemma 2b), 
it will be sufficient to consider the following two cases: 
\medskip 

a) $G=G^{\et }$; 
\newline 
In this case we can enlarge $O$ to assume that  
$G$ is a constant group scheme, where the property 
$G\simeq G[M]$ with $M\in\Mod _{O,\sigma }$ is obviously true.   
\medskip 

b) $G=G^{\loc }$.
\newline  
In this case enlarge $O$ to be able to assume that $G\otimes K$ is constant, 
i.e. $H=G(K_{\sep })$ has a trivial $\Gamma _K$-action. Choose a 
1-dimensional $\Bbb F_q$-subspace $H_2$ in $H$. It gives rise to 
a short exact sequence of $\Bbb F_q$-module group schemes 
$$0\longrightarrow G_2\longrightarrow G\longrightarrow G_1\longrightarrow 0$$
where $G_2(K_{\sep })=H_2$. 

Because $G$ is a strict $\Bbb F_q$-module, 
$\Bbb F_q$ acts on $t^*_G$ and, therefore, on $t^*_{G_2}$, 
via scalar multiplication. By Proposition 2,  
$G_2\simeq\mu _{\lambda }$, $\lambda\in O,\lambda\ne 0$, 
belongs to $\FGr '(\Bbb F_q)_O$ and, therefore, 
$G_1\in\FGr '(\Bbb F_q)_O$, cf. n.1.4. By induction, we have 
$G_1\simeq G[M_1]$, $M_1\in\Mod _{O,\sigma }$. 

Finally, $G$ satisfies the properties 1)-4) from n.2.5 and by Proposition 3, 
$G\simeq G[M]$ for $M\in\Mod _{O,\sigma }$. 
The Theorem is proved.
\medskip 

\ \

\subhead 3. Group schemes with strict $O_0$-action  
\endsubhead 
\medskip 

\ \ 

In this section $R=O_0$ is the valuation ring of a closed subfield 
$K_0$ of $K$ with finite residue field $\Bbb F_q$, $q=p^n$, 
$n\in\Bbb N$. Fix a choice of uniformising element $\pi _0$ in $K_0$. So, 
$O_0=\Bbb F_q[[\pi _0]]$.  

3.1. Denote by $\FGr '(O_0)_O$ a full subcategory  in $\FGr (O_0)_O$ 
consisting of objects killed by some power 
of $[\pi _0]$. This is an analogue of the classical category 
of $p$-torsion finite flat group schemes over the valuation ring 
of a  complete discrete valuation field of mixed characteristic. 

Let $\Mod (O_0)_{O}$ be the category, consisting of $M\in\Mod _{O,\sigma }$ 
with $O_0$-module structure given for $o\in O_0$ by endomorphisms 
$[o]\in\End _{\Mod _{O,\sigma }}(M)$ such that for all $m\in M$, 
$[\alpha ]m=\alpha m$ if $\alpha\in\Bbb F_q\subset k$ and $[\pi _0^N]m=0$ for 
sufficiently 
large $N$. Morphisms in the category 
$\Mod (O_0)_{O}$ are morphisms of the category 
$\Mod _{O,\sigma }$ which commute with all endomorphisms $[o]$, $o\in O_0$. 
In other words, any $M\in\Mod (O_0)_{O}$ is a free $O$-module $M$ 
of finite rank equiped with 
\medskip 

a) an $O$-linear map $\Phi :M_{\sigma }\longrightarrow M$ such 
that $\Phi\otimes K$ is isomorphism;
\medskip 

b) an $O$-linear nilpotent endomorphism $[\pi _0]:M\longrightarrow M$ 
commuting with $\Phi $. 
\medskip 

We shall denote by $\Mod '(O_0)_{O}$ a full subcategory of 
$\Mod (O_0)_{O}$ consisting of objects $M$ such that 
\medskip 

c) for any $m\in M$, $[\pi _0]m\equiv\pi _0m\operatorname{mod}\Phi (M_{\sigma })$.
\medskip 

The category $\Mod (O_0)_{O}$ is not generally abelian.  
But similarly to the category $\FGr '(O_0)_O$, one can work with short exact sequences 
in $\Mod (O_0)_{O}$. The sequence of objects and morphisms 
in $\Mod (O_0)_{O}$
$$0\longrightarrow M_1\mathbin{\overset{i}\to\longrightarrow} M
\mathbin{\overset{j}\to\longrightarrow }M_2\longrightarrow 0$$
is exact if and only if $i$ is a pure morphism of $O$-modules (i.e. $M/i(M_1)$ 
has no $O$-torsion) and $j$ is epimorphism of $O$-modules with 
kernel $i(M_1)$. This allows to define $O$-modules of 
equivalence classes of short exact sequences which satisfy the 
usual functorial properties.
\medskip 

3.2. Let $i:\FGr (O_0)_O\longrightarrow\FGr (\Bbb F_q)_O$ be the forgetful functor. 

\proclaim{Proposition 4} If $\Cal G\in\FGr '(O_0)_O$, then 
$i(\Cal G)\in\FGr '(\Bbb F_q)_O$.
\endproclaim

\demo{Proof} Indeed, if $\Cal G=(G,G^\flat )$ then $O_0$ acts 
via a scalar multiplication on $t^*_{G^\flat }$. This implies 
that $t^*_G$ is killed by some power of $\pi _0$ and, therefore, 
$t^*_{G\otimes K}=t^*_G\otimes K=0$, i.e. $G\otimes K$ is etale. 
The proposition is proved. 
\enddemo 

Let $L:\FGr '(\Bbb F_q)_O\longrightarrow\Mod _{O,\sigma }$ be the functor from n.2. 
For any $\Cal G\in\FGr '(O_0)_O$, $L(i(\Cal G))$ has a natural structure 
of an object of the category $\Mod (O_0)_{O}$. We denote by $L$ the 
functor from $\FGr '(O_0)_O$ to $\Mod (O_0)_{O}$ induced by the 
above correspondence $\Cal G\mapsto L(i(\Cal G))$. 

\proclaim{Theorem 2} The functor $L$ induces  antiequivalence of the categories 
$\FGr '(O_0)_O$ and $\Mod '(O_0)_{O}$.
\endproclaim 

\demo{Proof} Let $\Cal G\in\FGr '(O_0)_O$ and $M=L(\Cal G)\in\Mod (O_0)_{O}$. 
Considering $M$ as an object of $\Mod _{O,\sigma }$ we can recover the 
$\Bbb F_q$-module scheme $i(\Cal G)$ in the form $G[M]$. This reduces the proof 
of our theorem to the following statement:
\medskip 
\proclaim{Proposition 5} $\Cal G$ is a strict $O_0$-module iff $M$ satisfies 
the condition {\rm c) from n.3.1.}
\endproclaim 
\medskip 

Let $\Cal G^{loc}$ be the maximal local subobject in $\Cal G$ 
and $M^{loc}=L(\Cal G^{loc})$. 
Then the natural projection $M\longrightarrow M^{loc}$ induces isomorphism 
of $(O_0,O)$-modules $M/\Phi (M_{\sigma })$ and 
$M^{loc}/\Phi (M_{\sigma }^{loc})$. So, it will be sufficient to prove 
the above proposition for a local $\Cal G$. Under this assumption 
choose a vector column  $\bar m=(m_1,\dots ,m_n)^t$ such that its coordinates create 
an $O$-basis of  $M=L(\Cal G)$ and set $\Phi \bar m=C\bar m$ and 
$[\pi _0]\bar m=D\bar m$, where $C,D\in\Bbb M_n(O)$ are $n\times n$-matrices 
with coefficients from $O$.

\proclaim{Lemma 4} The above matrices $C$ and $D$ give a structure of an object of 
the category $\Mod '(O_0)_{O}$ on $M$ if and only if 

{\rm 1)} $\operatorname{det}C\ne 0$;

{\rm 2)} $D$ is nilpotent and $\sigma (D)C=CD$;

{\rm 3)} there is a $B\in\Bbb M_n(O)$ such that $D-\pi _0E_n=BC$.
\endproclaim 

\demo{Proof} Indeed, 1) means that $\Phi (M_{\sigma })$ is 
a lattice in the $O$-module $M$.

Property 2) means that $M$ is $[\pi _0]$-torsion and the morphisms 
$\Phi $ and $[\pi _0]$ commute one with another. Indeed, 
$(\Phi\circ[\pi _0])(\bar m)=[\pi _0](\Phi\bar m)=CD\bar m$ and 
$([\pi ]\circ\Phi )(\bar m)=\sigma (D)C\bar m$. 

Notice that the coordinates of the vector $C\bar m$ generate 
the $O$-module $\Phi (M_{\sigma })$. So, the coordinates 
of $[\pi _0](\bar m)-\pi _0\bar m$ belong to $\Phi (M_{\sigma })$ 
if and only if there is an $B\in\Bbb M_n(O)$ such that 
$D-\pi _0E_n=BC$. 

The lemma is proved.
\enddemo 

\remark{Remark} a) $M=M^{loc}$ if and only if $C$ is $\sigma $-nilpotent, 
i.e. there is an $N\in\Bbb N$ such that $\sigma ^{N}(C)\dots\sigma (C)C=0$;
\newline 
b) the above properties 1) and 3) imply that $\sigma (D)-\pi _0E_n=CB$.
\endremark 
\medskip 

With the above notation and agreements let 
$\Cal G=(G,G^\flat )$, where $A(G)=O[\bar X]/I$ and $A(G^\flat )=O[\bar X]/(I\cdot I_0)$, 
where $\bar X=(X_1,\dots ,X_n)$ and the ideal $I$ is generated by the coordinates of 
the vector $\bar X^q-C\bar X$. 

Then $A(G)$ has an $O$-basis $\{X_1^{a_1}\dots X_n^{a_n}\ |\ 0\leqslant a_1,\dots ,a_n<q\}$ 
and this basis can be completed to an $O$-basis of $A(G^\flat )$ by 
joining the elements of the set 
\linebreak 
$\{X_i^q\ |\ 1\leqslant i\leqslant n\}$. 

The action of $O_0$ is given by $[\pi _0]\bar X=D\bar X$ on $A(G)$. 
This action is strict if and only if it can be extended to $A(G^\flat )$ 
via the relation 
$$[\pi _0]\bar X=D\bar X+B(\bar X^q-C\bar X)$$
with some $B\in\Bbb M_n(O)$ in such a way that it induces 
the scalar multiplication by $\pi _0$ 
on $t^*_{A(G^\flat )}$ and $N_{A(G^\flat )}$. 

The first condition is equivalent to the equality $D-BC=\pi _0E_n$.

The second can be analysed as follows,  
$$[\pi _0]^\flat (\bar X^q-C\bar X)\equiv  
([\pi _0]^\flat\bar X)^q-C([\pi _0]^\flat\bar X)
\equiv (\sigma (D)-CB)(\bar X^q-C\bar X)$$ 
and is equivalent to the matrix equality $\sigma D-CB=\pi _0E_n$. 

So, Proposition 5 and Theorem 2 follow from Lemma 4 and the above Remark b).  
\enddemo 

3.3. Clearly, the antiequivalence $L$ transforms short exact sequences in $\FGr '(O_0)_O$ to 
short exact sequences in the category $\Mod '(O_0)_{O,\sigma }$. 
In particular, we have the following property.

\proclaim{Proposition 6} Suppose $\Cal G_1, \Cal G_2\in\FGr '(O_0)_O$ and 
$f\in\Hom _{\FGr'(O_0)_{O}}(\Cal G_1,\Cal G_2)$. Then 
\newline 
{\rm a)} $f$ is a closed embedding if and only if $L(f)$ is surjective;
\newline 
{\rm b)} $f$ is fully faithful if and only if $L(f)$ is a pure embedding. 
\endproclaim

Suppose for $n\in\Bbb N$, $\Cal G^{(n)}\in\FGr '(O_0)_O$, 
$i_n:\Cal G^{(n)}\longrightarrow\Cal G^{(n+1)}$ is a closed 
immersion, $j_n:\Cal G^{(n+1)}\longrightarrow\Cal G^{(n)}$ is a 
fully faithful morphism. For $m>n$, set 
\linebreak 
$i_{nm}=i_{n}\circ\dots \circ i_{m-1}$ and 
$j_{mn}=j_{m-1}\circ\dots\circ j_n$. Then 
(following the original deinition of Tate) 
$\{\Cal G^{(n)}, i_n,j_n\}$ is 
a $\pi _0$-divisible 
group in the category $\FGr '(O_0)_O$ if for any $m>n$, 
$(\Cal G^{(n)},i_{nm})=\Ker [\pi _0^m]\id _{\Cal G^{(m)}}$ 
and $[\pi _0^{m-n}]\id _{\Cal G^{(m)}}=j_{mn}\circ i_{nm}$. 

Similar definitions can be done for the category $\Mod '(O_0)_{O}$, 
where a $\pi _0$-divisible group appears as the collection 
$\{M^{(n)}, i_n,j_n\}$ where $i_n$ is a pure embedding of underlying $O$-modules 
$M^{(n)}\longrightarrow M^{(n+1)}$ and $j_n:M^{(n+1)}\longrightarrow M^{(n)}$ 
is a surjection of 
$O$-modules. 

3.4. Clearly, the functor $L$ from Theorem 2 transforms $\pi _0$-divisible groups 
in the category $\FGr '(O_0)_O$ to $\pi _0$-divisible groups in 
the category $\Mod '(O_0)_{O}$. 

\proclaim{Theorem 3}  For any object $\Cal G$ of the category $\FGr '(O_0)_O$, 
there is a $\pi _0$-divisible 
group $\{\Cal H^{(n)}, i_n,j_n\}_{n\geqslant 1}$  and a closed embedding 
$\Cal G\longrightarrow\Cal H^{(N)}$, where $N\in\Bbb N$ is such that 
$[\pi _0^N]\id _{\Cal G}=0$. 
\endproclaim 

\remark{Remark} The statement of the above theorem is equivalent 
to the existence of a $\pi _0$-divisible 
group $\{{\Cal H'}^{(n)},i_n,j_n\}_{n\geqslant 1}$ and a fully faithful morphism 
${\Cal H'}^{(N)}\longrightarrow\Cal G$. 
\endremark 

\demo{Proof} The antiequivalence $L$ allows us to prove the dual version of 
this theorem in 
the category $\Mod '(O_0)_{O}$. 

Let $M\in\Mod '(O_0)_{O}$ and let $N\in\Bbb N$ be such that 
$[\pi _0^N]M=0$. Use induction on $N$. 

Suppose first, that $N=1$.

Then $M$ is a free $O$-module with $\sigma $-linear 
$\Phi :M_{\sigma }\longrightarrow M$ such that 
$\pi _0M\subset\Phi (M_{\sigma })\subset M$. Choose an 
$O$-basis $m_1,\dots ,m_n$ in $M$ and take a vector column 
$\bar m=(m_1,\dots ,m_n)^t$. 
Then $\Phi\bar m=C\bar m$, where $C\in\Bbb M_n(O)$ is a divisor of 
$\pi _0E_n$. Let $\widetilde{C}\in\Bbb M_n(O)$ be such that 
$C\widetilde{C}=\pi _0E_n$. 

For $l\in\Bbb N$, introduce free $O$-modules $M^{(l)}$ with free generators 
$$\{m^{(k)}_{1i}, m^{(k)}_{2i}\ |\ 1\leqslant i\leqslant n, 1\leqslant k\leqslant l\}$$

Let $\bar m_1^{(0)}=\bar m_2^{(0)}=\bar 0$ and for $k\geqslant 1$, let 
$\bar m_1^{(k)}=(m_{11}^{(k)},\dots ,m_{1n}^{(k)})^t$ and 
\linebreak  
$\bar m_2^{(k)}=(m_{21}^{(k)},\dots ,m_{2n}^{(k)})^t$. Define $O$-linear 
morphisms $\Phi :M_{\sigma }^{(l)}\longrightarrow M^{(l)}$ 
and $[\pi _0]:M\longrightarrow M$ by setting 
for $1\leqslant k\leqslant l$, 

$$\Phi \bar m_1^{(k)}=\widetilde{C}\bar m_1^{(k)}+\bar m_2^{(k-1)},\ \ 
\Phi \bar m_2^{(k)}=C\bar m_2^{(k)}+\bar m_1^{(k)} $$
and $[\pi _o]\bar m_1^{(k)}=\bar m_1^{(k-1)}$, $[\pi _0]\bar m_2^{(k)}=\bar m_2^{(k-1)}$. 

It is easy to see that we defined a structure of objects of the category 
$\Mod '(O_0)_{O}$ on all $M^{(l)}$, $l\in\Bbb N$, and 
the system $\{M^{(l)}\}_{l\geqslant 1}$ together with the natural 
inclusions $i_n:M^{(n}\longrightarrow M^{(n+1)}$ and projections 
$j_l:M^{(l+1)}\longrightarrow M^{(l)}$ gives a $\pi _0$-divisible group in 
the category $\Mod '(O_0)_{O}$. Clearly, the correspondences 
$\bar m_1^{(1)}\mapsto \bar 0$ and $\bar m_2^{(1)}\mapsto \bar m$ give 
an epimorphic map from $M^{(1)}$ to $M$. The case $N=1$ has been  
considered.

Suppose $N>1$ and Theorem 3 has been proved for all $M'\in\Mod '(O_0)_{O}$ 
such that $[\pi _0]^{N-1}(M')=0$. 

Let $M_1=\Ker [\pi _0]^{N-1}\id _M$ and $M_2=[\pi _0]^{N-1}(M)\subset M$. 
Then $M_1$ and $M_2$ have natural structures of objects of the category 
$\Mod '(O_0)_{O}$ and in this category we have a natural short exact sequence 
$$\varepsilon :\ \ 0\longrightarrow M_1\longrightarrow M\longrightarrow M_2\longrightarrow 0$$
(Notice that, generally, the embedding $M_2\subset M$ is not pure). 
By induction, there is a $\pi _0$-divisible group 
$\{T^{(n)},i_n,j_n\}_{n\geqslant 1}$ and pure embeddings 
$\alpha :M_1\longrightarrow T^{(N-1)}$, 
$\beta :M_2\longrightarrow T^{(1)}$. Consider the short exact sequences 
$$\alpha _*\varepsilon :\ \ \ 
0\longrightarrow T^{(N-1)}\longrightarrow\alpha _*M
\mathbin{\overset{j_\alpha }\to\longrightarrow} 
M_2\longrightarrow 0$$
and
$$\beta ^*\eta _N :\ \ \ 
0\longrightarrow T^{(N-1)}\longrightarrow\beta ^*T^{(N)}
\mathbin{\overset{j_\beta }\to\longrightarrow}M_2\longrightarrow 0$$
where the second sequence  is obtained via $\beta $ from the standard short exact sequence 
$$\eta _N:\ \ 0\longrightarrow T^{(N-1)}\mathbin{\overset{i_N}\to\longrightarrow}
T^{(N)}\mathbin{\overset{j_{N1}}\to\longrightarrow} T^{(1)}\longrightarrow 0$$
Notice that there is a pure embedding $M\longrightarrow \alpha _*M$. Notice also that 
in the both sequences 
$\alpha _*\varepsilon $ and $\beta ^*\eta _N$, the epimorphic maps $j_{\alpha }$ and 
$j_{\beta }$ are 
induced by multiplication by $[\pi _0^{N-1}]$. In other words, 
the endomorphism $[\pi _0^{N-1}]$ on $\alpha _*M$ and, resp., on  $\beta ^*T^{(N)}$ 
is a  composition of 
$j_\alpha $ and, resp.,  of $j_{\beta }$ with the natural inclusion of $M_2$ into 
$\alpha _*M$ and, resp., into  $\beta ^*T^{(N)}$ 
(which is  induced by the embedding $M_2\subset T^{(1)}$). 

Let 
$$0\longrightarrow T^{(N-1)}\longrightarrow\widetilde{M}\longrightarrow M_2\longrightarrow 0$$ 
be the difference of $\alpha _*\varepsilon $ and $\beta ^*\eta _N$ 
considered as elements of the group 
\linebreak  
$\Ext _{\Mod '(O_0)_{O}}(M_2,T^{(N-1)})$. It is easy to see 
that $[\pi _0^{N-1}](\widetilde M)=0$ and 
there is an $\widetilde{M}_1\in\Mod '(O_0)_{O}$, a surjective 
map $\tilde\jmath$ and a pure embedding $\tilde\imath$ such that  
$$\beta ^*T^{(N)}\oplus\widetilde{M}
\mathbin{\overset{\tilde\jmath}\to\longrightarrow}\widetilde {M}_1
\mathbin{\overset{\tilde\imath}\to\longleftarrow}\alpha _*(M)$$

By inductive assumption,  there is 
a pure embedding of $\beta ^*T^{(N)}\oplus\widetilde{M}$
into a $\pi _0$-divisible group in the category 
$\Mod '(O_0)_{O}$. This implies easily the existence of such 
an embedding for $\alpha _*M$ and, therefore, for $M$.

Theorem 3 is proved. 
\enddemo 
\medskip 

3.5. {\it Relation to the mixed characteristic case.}
\medskip 

In nn.3.5.1-3.5.3 below we need a full subcategory 
$\FGr '_1(O_0)_O$ in $\FGr '(O_))_O$. It consists of 
$\Cal G=(G,G^\flat )$ such that $[\pi _0]\id _G=0$. 
Then the functor $L$ induces an antiequivalence of this category 
and the full subcategory $\Mod '_1(O_0)_O$ of $\Mod '(O_0)_O$ 
consisting of free $O$-modules $M$ such that 
$\pi _0M\subset\Phi (M_{\sigma })\subset M$. As usually, 
$e=e(K/K_0)$ is the ramification index of $K=\Frac O$ over $K_0=\Frac O_0$.
\medskip 

3.5.1. Suppose $e(K/K_0)=1$. 

Introduce the category $\SH _1(\Bbb F_q)_O$ with objects 
$(M^0,M^1,\varphi _0,\varphi _1)$, where 
$M^0$ is an $O$-module of finite length such that $\pi _0M=0$, 
$\varphi _0:M^0\longrightarrow M^0$ is a $\sigma $-linear morphism, 
$M^1=\Ker \varphi _0$, $\varphi _1:M^1\longrightarrow M^0$ is a 
$\sigma $-linear morphism and $\varphi _0(M^0)+\varphi _1(M^1)=M$. 
This is an analogue of Fontaine's category of 
filtered modules. 

Consider the functor $\SH :\Mod _1'(O_0)_O\longrightarrow\SH _1(\Bbb F_q)$ 
defined by the correspondence $M\mapsto (M^0,M^1,\varphi _0,\varphi _1)$, 
where $M^0=M\operatorname{mod}\pi _0M$, $\varphi _0=\Phi\operatorname{mod}\pi _0$ 
and $\varphi _1$ is induced by $\dsize\frac{1}{\pi _0}\Phi $. The functor 
$\SH $ is an equivalence of categories if $q>2$ and is ``very close'' to an 
equivalence if $q=2$. This shows that strict $O_0$-modules have similar description 
as conventional group schemes if $e=1$. Actually, one can develop 
the Dieudonne theory in the context of strict modules and realise 
the approach from [Ab1] to recover directly an analogue of Fontaine's classification 
of group schemes over Witt vectors. 
\medskip 

3.5.2. Suppose $e\leqslant q-1$. 

Define the category $\SH _1(\Bbb F_q)_O$ of the collections 
$(M,M^0,M^1,\varphi _0,\varphi _1)$, where $M$ is an $O$-module 
of finite length killed by $\pi _0$, $M^0=\Ker \pi |_M$, $\varphi _0:M^0\longrightarrow M$ 
is a $\sigma $-linear map, $M^1=\Ker \varphi _0$, $\varphi _1:M^1\longrightarrow M$ 
is a $\sigma $-linear map, and $\varphi _0(M^0)\otimes _kO+\varphi _1(M^1)\otimes _kO=M$.
This category is an analogue of the category $\SH _O$ from [Ab3]. Consider the functor 
$\SH :\Mod _1'(O_0)_O\longrightarrow\SH (\Bbb F_q)_O$ defined by the correspondence 
$M\mapsto (\bar M,\bar M^0,\bar M^1,\varphi _0,\varphi _1)$, where 
$\bar M$ is the $O$-submodule in 
$\dsize\frac{\pi }{\pi _0}O\otimes M\operatorname{mod}\pi M$ 
generated by the images of the elements of the sets 
\linebreak 
$\dsize\{\frac{\pi }{\pi _0}\Phi (m\ |\  m\in M\}$ and 
$\dsize\{\frac{1}{\pi _0}\Phi (m)\ |\  m\in M, \ \Phi (m)=0\}$. 
Then $\bar M^0=
\Ker\pi |_{\bar M}=M\operatorname{mod}\pi M$, $\varphi _0:
\bar M^0\longrightarrow M$ is induced by 
$\dsize\frac{\pi }{\pi _0}\Phi $, $\bar M^1=\Ker\varphi _0$ and 
$\varphi _1:\bar M^1\longrightarrow M$ is induced by 
$\dsize\frac{1}{\pi _0}\Phi $. 

If $e<q-1$ then $\SH $ is an equivalence of categories and, if $e=q-1$ 
then $\SH $ is ``very close'' to an equivalence of categories. 
Again, the methods from [Ab3] 
can be used to obtain directly the classification of objects from 
$\FGr '_1(O_0)_O$ via objects of the category 
$\SH (\Bbb F_q)_O$. (The details will appear in W.Gibbons's thesis.)

Notice that when working with equal characteristic case the category $\SH (\Bbb F_q)_O$ 
can be replaced by a simpler category 
$\SH '(\Bbb F_q)_O$ consisting of the triples 
$(\bar M,\bar M^0,\varphi )$, where $\bar M$ is an $O$-module of finite length killed 
by $\pi _0$, $\bar M^0=\Ker \pi |_{\bar M}$ and 
$\varphi _0:\bar M^0\longrightarrow \bar M$ 
is a $\sigma $-linear morphism such that 
$\varphi (\bar M^0)\otimes O=\bar M$. The functor $\SH ':\FGr '_1(O_0)_O
\longrightarrow\SH '(\Bbb F_q)_O$ is defined by the correspondence 
$M\mapsto (\bar M,\bar M^0,\varphi _0)$, where 
$\bar M$ is an $O$-submodule in $\dsize\frac{1}{\pi _0}O\otimes M
\operatorname{mod}\pi M$ generated by elements of the form 
$\dsize\frac{1}{\pi _0}\Phi m$ with $m\in M$, 
$\bar M^0=\Ker \pi |_M=M\operatorname{mod}\pi M$ and $\varphi _0$ is induced by 
$\dsize\frac{1}{\pi _0}\Phi $ (this map is still $\sigma $-linear because of 
the equal characteristic situation). 
\medskip

3.5.3. Suppose that $e(K/K_0)$ is arbitrary. 

Introduce the category $\BR _1(\Bbb F_q)_O$. Its objects are 
the triples $(M,M^1,\varphi )$, where $M$ is 
an $O$-module of finite length such that $\pi _0M=0$, 
$M^1$ is an $O_1$-submodule in $M\otimes _OO_1$, where $O_1=O[\pi _1]$ with $\pi _1^p=\pi $, 
and $\varphi :M^1\longrightarrow M$ is a $\sigma $-linear map 
($\sigma $ is still the $q$-th power map) such that 
$\varphi (M^1)=M$. This category is an equicharacteristic version 
of Breuil's category $\Mod _{/S_1}$ appeared in his classification of 
period $p$ group schemes in the mixed characteristic case, cf. [Br].  

Again there is a natural functor from 
$\Mod _1'(O_0)_O$ to $\BR _1(\Bbb F_q)_O$ defined by the correspondence 
$M\mapsto (\bar M,\bar M^1,\varphi )$, where 
$\bar M=M\operatorname{mod}\pi _0$ and $\varphi $ is induced 
by $\dsize\frac{1}{\pi _0}\Phi $. But one can't expect that Breuil's method 
works in the equal characteristic case, because it is based 
very heavily on crystalline technique. 
\medskip 

3.5.4. Again $e(K/K_0)$ is arbitrary.

Introduce an equal characteristic analogue of the concept 
of $p$-etale $\varphi $-module of $q$-heighjt 1 over $S$ from [Fo4]. 

First, introduce 2nd copy $\underline{O}_0=\Bbb F_q[[\underline{\pi }_0]]$ of 
$O_0$ (actually, $\underline{\pi }_0$ is an analogue of $[\pi _0]$). 
Then $\underline{O}_0$ will be considered as an analogue of 
$\Bbb Z_p$ and $\underline{\widetilde{O}}_0=k[[\underline{\pi }_0]]$ 
is an analogue of 
Witt vectors $W(k)$. Then consider 
$S=\underline{\widetilde{O}}_0[[\pi ]]$, where $k[[\pi ]]=O$ appears 
as a subring of $S$. Introduce $\sigma :S\longrightarrow S$ such that 
$\sigma |_{\underline{O}_0}=\id $ and $\sigma |_O$ is 
(as earlier) the $q$th power map.

Suppose that $M$ is an $S$-module of finite type together with an 
$\underline{O}_0$-linear map $\Phi :M_{\sigma }:=
M\otimes _{(O,\sigma )}O\longrightarrow M$. Then $M$ is called 
a $\underline{\pi} _0$-etale $\varphi $-module of $q$-height 
1 (with $q=\underline{\pi} _0-\pi _0$), if $\operatorname{Coker}\Phi $ is 
killed by multiplication by $\underline{\pi} _0-\pi _0$.

With this notation Theorem 2 establishes a description of the category 
of finite flat strict $O_0$-modules over $O$ in terms of $\underline{\pi} _0$-torsion 
$\underline{\pi} _0$-etale $\varphi $-modules of $q$-height 1. 
\medskip

\subhead 4. Properties of arising Galois modules 
\endsubhead 
\medskip 

\ \ 

In this section $\Cal G=(G,G^\flat )\in\FGr '(O_0)_O$, 
$H=G(K_{\sep })$ is $O_0[\Gamma _K]$-module of geometric points 
of $G$ and $e=e(K/K_0)$ --- is the ramification index of 
$K$ over $K_0$. We also set 
$\Gamma _K=\Gal (K_{\sep }/K)$ and denote by $I_K$ the inertia 
subgroup of $\Gamma _K$. 
\medskip 

\subsubhead {\rm 4.1.} Characters of the semisimple envelope of $H$ 
\endsubsubhead

Suppose $\bar k$ is an algebraic closure of $k$ and the character 
$\chi :I_K\longrightarrow\bar k^*$ appears with a nonzero multiplicity 
in the semisimple envelope 
of the $O_0[\Gamma _K]$-module $H$. An analogue of the Serre Conjecture for $H$ 
can be stated as follows.

\proclaim{Theorem 4} For the above 
character $\chi $, there are $a,N\in\Bbb N\setminus p\Bbb N$ 
such that $\chi =\chi _N^a$, where $a=a_0+a_1q+\dots +a_{N-1}q^{N-1}$ 
with $0\leqslant a_i\leqslant e$ and $\chi _N:I_K\longrightarrow\bar k^*$ 
is such that for any $\tau\in I_K$, $\chi _N(\tau )=\tau (\pi _N)/\pi _N$,  
where $\pi _N\in K_{\sep }$ and $\pi _N^{q^N-1}=\pi $. 
\endproclaim 

\demo{Proof} This can be deduced in the same way as it has been 
obtained in the case of usual group schemes in [Ra]. 
First we can assume that $k=\bar k$ and $e<q-1$. 
Then any simple object of 
the category $\Mod '(O_0)_O$ appears in the form 
$M=\oplus _{0\leqslant i<n}Om_i$, where  
$\Phi m_0=\pi ^{a_0}m_1,\dots ,\Phi m_{N-1}=\pi ^{a_{N-1}}m_0$ 
and $[\pi _0]m_0=\dots =[\pi _0]m_{N-1}=0$ with 
$0\leqslant a_i\leqslant e$, $0\leqslant i<N$.  

Then the corresponding Galois module consists of $K_{\sep}$-points 
of the $O$-algebra 
$O[T_0,\dots ,T_{N-1}]$, where 
$T_0^q=\pi ^{a_0}T_1,\dots ,T_{N-1}^q=\pi ^{a_{N-1}}T_0$. It can be 
naturally identified with the 
$O_0[\Gamma _K]$-module $\{\alpha\pi _N^a\ |\ \alpha\in\Bbb F_{q^N}\}$, 
where $a=a_0+a_1q+\dots +a_{N-1}q^{N-1}$.  
Clearly, $I_K$ acts on it via the conjugacy class of 
characters 
\linebreak 
$\{\sigma ^i\chi _N^a\ |\ 0\leqslant i<N\}$. 
\enddemo 

\remark{Remark} Following Raynaud's method one can deduce 
from the above description of simple objects in 
$\Mod '(O_0)_O$ that if $e<q-1$ then the functor 
$G\mapsto G(K_{\sep})$ is a fully faithful 
functor from $\Mod '(O_))_O$ to the category of 
$O_0[\Gamma _K]$-modules. 
\endremark 
\medskip 

\subsubhead {\rm 4.2.} Ramification estimates 
\endsubsubhead 
\medskip 

These estimates are given in Theorem 5 below and are completely 
similar to the known estimates in the case of conventional 
group schemes, cf. [Fo2]. The proof is based on the knowledge of 
\lq\lq equations\rq\rq\ of the strict module $G$ and is done below by 
the methods of the paper [Ab4]. Notice that the methods from 
[Fo3] also can be adjusted to obtain the same estimates.

\proclaim{Theorem 5} If $H$ is killed by $[\pi _0^N]$ then 
the ramification subgroups $\Gamma _K^{(v)}$ 
act trivially on $H$ for $v>e\dsize \left (N+\frac{1}{q-1}\right )-1$.
\endproclaim 

\demo{Proof} We can assume that there is a $\pi _0$-divisible group 
$\{\Cal G^{(i)}\}_{i\geqslant 1}$ of a height $h$ in 
$\FGr '(O_0)_O$ such that $\Cal G^{(N)}=\Cal G$. 

4.2.1. By the above assumption,  
$L(\Cal G)$ is a free $O$-module of rank $hN$ and we can choose 
its $O$-basis in the form 
$$m_1,\dots ,m_h,[\pi _0]m_1,\dots ,[\pi _0]m_h,\dots ,
[\pi _0]^{N-1}m_1,\dots ,[\pi _0^{N-1}]m_h$$

Introduce vector-columns 
$$\bar m_1=(m_1,\dots ,m_h)^t,\dots ,
\bar m_N=[\pi _0^{N-1}]\bar m_1=([\pi _0^{N-1}]m_1,\dots ,[\pi _0^{N-1}]m_h)^t$$

Then, in the obvious notation, 
the structure of $M=L(\Cal G)$ can be given  
in the following form

$$\frac{1}{\pi _0}C_1\Phi\bar m_N+\frac{1}{\pi _0}C_2\Phi\bar m_{N-1}+
\dots +\frac{1}{\pi _0}C_{N}\Phi\bar m_{1}=
\bar m_N-\frac{1}{\pi _0}\bar m_{N-1}$$ 
$$\qquad\qquad\frac{1}{\pi _0}C_1\Phi \bar m_{N-1}+\dots +\frac{1}{\pi _0}C_{N-1}
\Phi\bar m_1=\bar m_{N-1}-\frac{1}{\pi _0}\bar m_{N-2}$$
$$ ...............................$$
$$\qquad\qquad\qquad\qquad\qquad\frac{1}{\pi _0}C_1\Phi\bar m_1=\bar m_1$$
where all $C_i\in\GL (h,O)$ and $\det C_1\ne 0$. 

Consider vector columns $\bar X_i=(X_{i1},\dots ,X_{ih})^t$ of 
independent variables $X_{ij}$, 
\linebreak 
$1\leqslant i\leqslant N$, 
$1\leqslant j\leqslant h$. Then the algebra $A(G)$ appears as 
a quotient of $O[\bar X_1,\dots ,\bar X_N]$ by the ideal generated 
by equations 
$$\sum\Sb 1\leqslant i\leqslant s\endSb C_i\bar X_{s+1-i}^q=\pi _0\bar X_s-\bar X_{s-1}\tag{5}$$
where $1\leqslant s\leqslant N$ and by definition $\bar X_0=\bar 0$. 

Consider the points of $G(K_{\sep})$ as solutions 
$\bar a=(\bar a_1,\dots ,\bar a_N)$ of the system (5). 
The following lemma can be easily proved by induction on $N$.

\proclaim{Lemma 5} If $\bar a=(\bar a_1,\dots ,\bar a_N)$ and 
$\bar a'=(\bar a_1',\dots ,\bar a_N')$ are solutions of $(4.1)$ 
such that $\bar a\equiv\bar a'\operatorname{mod}\pi _0^{\frac{1}{q-1}}m_{\sep}$, 
where $m_{\sep}$ is the maximal ideal of the valuation ring of $K_{\sep}$, 
then $\bar a=\bar a'$.
\endproclaim 
\medskip 

4.2.2. Suppose $\alpha\in\Bbb Q_{>0}$ has the $p$-adic valuation 0. 
Then $\alpha =\dsize\frac{m}{q^M-1}$ with suitable $m,M\in\Bbb N$, $(m,p)=1$. 
For any such $\alpha $ there is an extension $K_{\alpha }$ of $K$ with 
$[K_{\alpha }:K]=q^M$ and the Herbrand function 

$$\varphi _{K_{\alpha }/K}(x)=\cases x, &\text {for $0\leqslant x\leqslant \alpha $} \\ 
\alpha +\dsize\frac{x-\alpha}{q^M}, & \text {for $x\geqslant\alpha $}. \endcases 
$$
Notice that $\varphi _{K_{\alpha }/K}$ has only one edge point in $x=\alpha $. 

Explicit construction of $K_{\alpha }$ can be found in [Ab4], where the field $K_{\alpha }$ 
appears as an extension of $K$ of degree $q^M$ in $L_{\alpha }=K(\pi _M)(T)$, 
where $\pi _M^{q^M-1}=\pi $ and $T^{q^M}-T=\pi _M^{-m}$. Clearly, $K_{\alpha }$ 
is totaly ramified over $K$ and, therefore, there is a field isomorphism 
$$h_{\alpha }:K\longrightarrow K_{\alpha }$$
From the above construction of $K_{\alpha }$, it follows easily that $h_{\alpha }$ can be 
chosen in such a way that for any $a\in m_K$ 
($m_K$ is the maximal ideal in $O$), 
$$a=h_{\alpha }(a)^{q^M}+\tilde a,$$
with $\tilde a\in K_{\alpha }$ such that  $v_K(\tilde a)\geqslant v_K(a)+\alpha $ 
($v_K$ is the normalized valuation in $K$). 
\medskip 

4.2.3. Denote by the same symbol an extension of $h_{\alpha }$ to an isomorphism 
of $K_{\sep }$ onto $K_{\alpha ,\sep}=K_{\sep}$. Clearly, 
$\bar X=(\bar X_1,\dots ,\bar X_s)\mapsto 
h_{\alpha }(\bar X)=
(h_{\alpha }(\bar X_1),\dots ,h_{\alpha }(\bar X_s))$ is a one-one correspondence between 
solutions of the system (5) and solutions 
$\bar Y=(\bar Y_1,\dots ,\bar Y_s)$ of the similar system

$$\sum\Sb 1\leqslant i\leqslant s\endSb h_{\alpha }(C_i)\bar Y_{s+1-i}^q=
h_{\alpha }(\pi _0)\bar Y_s-\bar Y_{s-1}\tag{6}$$
where $1\leqslant s\leqslant N$ and by definition $\bar Y_0=\bar 0$. 

\proclaim{Lemma 6} If $\alpha >\dsize e\left (N+\frac{1}{q-1}\right )-1$, then 
for any solution $\bar X^{(0)}$ of (5) there is a unique solution 
$\bar Y^{(0)}$ of (6) such that 
$\bar X^{(0)}\equiv \bar Y^{(0)q^M}\operatorname{mod}\pi _0^{\frac{1}{q-1}}m_{\sep }$
\endproclaim  

\demo{Proof of lemma} 
the correspondence $\bar Y\mapsto\bar Z=\bar Y^{q^M}-\bar X^{(0)}$ establishes a 
one-one correspondence between solutions $\bar Y$ of (5) and solutions 
$\bar Z=(\bar Z_1,\dots ,\bar Z_s)$ of the system of equations 
$$\sum\Sb 1\leqslant i\leqslant s\endSb C_i\bar Z_{s+1-i}^q=
\pi _0\bar Z_s-\bar Z_{s-1}+\bar F_s\tag{7}$$
where $1\leqslant s\leqslant N$, $\bar Z_0=0$ and 
$$\bar F_s=\tilde\pi _0\bar Y_s^{q^M}-\sum\Sb 1\leqslant i\leqslant s\endSb 
\tilde C_i\bar Y_{s+1-i}^{q^{M+1}}\in\pi _0^{1+\frac{1}{q-1}}m_{\sep }$$
because $v_K(\tilde \pi _0),v_K(\tilde C_i)\geqslant 1+
\alpha>e\dsize\left (1+\frac{1}{q-1}\right )$. 

Now induction on $s$ shows that the system (7) has a unique 
solution $\bar Z$ with coordinates in $\pi _0^{\frac{1}{q-1}}m_{\sep }$.

Lemma is proved.
\enddemo 

With the above notation and assumptions we have the following corollary.

\proclaim{Corollary} Suppose $E$, resp. $E_{\alpha }$, is obtained by joining to 
$K$, resp. $K_{\alpha }$, all coordinates of all solutions of the system of 
equations (4.1), resp. (4.2), in $K_{\sep }$. Then $EK_{\alpha }=E_{\alpha }$. 
\endproclaim 
\medskip 

4.2.4. For any finite extension $F\subset L$ in $K_{\sep }$, let 
$v(L/F)$ be the minimal rational number such that 
the ramification groups 
$\Gamma _F^{(v)}$ act trivially on $L$ for $v>v(L/F)$. 

\proclaim{Proposition 7}  With the above notation 
there is the following inequality 
$$v(E/K)\leqslant e\left (N+\frac{1}{q-1}\right )-1$$
\endproclaim 

\demo{Proof} 
Suppose this inequality does not hold. Then there is a rational number $\alpha $ 
satisfying the assumptions from the beginning of n.4.2.2 
and the inequalities  
$$v(E/K)>\alpha >e\left (N+\frac{1}{q-1}\right )-1.$$

Notice that  $E_{\alpha }=EK_{\alpha }$ implies that  
$$v(E_{\alpha }/K)=\max\{v(E/K),v(K_{\alpha }/K)\}=v(E/K)$$

On the other hand, looking at the 
maximal edge points of Herbrand functions from the identity 
$\varphi _{E_{\alpha }/K}=\varphi _{E_{\alpha }/K_{\alpha }}\circ\varphi _{K_{\alpha }/K}$, 
we obtain that 

$$v(E_{\alpha }/K)=\max\left\{ v(K_{\alpha }/K),
\varphi _{K_{\alpha }/K}(v(E_{\alpha }/K_{\alpha }))\right\}=$$ 
$$\max\left\{\alpha ,\frac{v(E/K)-\alpha }{q}+\alpha\right\}<v(E/K)$$
because $v(E/K)=v(E_{\alpha }/K_{\alpha })$. Contradiction. 

Theorem 5 is proved.
\enddemo 
\medskip 

4.3. As it was noticed in n.3.5.1, if $e=1$ 
then killed by $[\pi _0]$ strict $O_0$-modules  
behave very similarly to group schemes of period $p$ 
over Witt vectors. For this reason, one can apply 
directly methods from [Ab2] to prove the following result.

\proclaim{Theorem 6} 
Suppose $H$ is an $\Bbb F_q[\Gamma _K]$-module such that 
\newline 
{\rm a)} the action of inertia subgroup of 
$\Gamma _K$ on the semisimple envelope of $H$ is given by 
characters, which satisfy Serre's Conjecture, cf. Theorem 4;
\newline 
{\rm b)} the ramification subgroups $\Gamma _K^{(v)}$ act trivially 
on $H$ if $v>\dsize\frac{1}{q-1}$ (i.e. the ramification estimate 
from Theorem 5 holds for $H$).

Then there is an $\Cal G=(G,G^\flat )\in\FGr '(O_0)_O$ 
such that $H\simeq G(K_{\sep})$.
\endproclaim 
\enddemo

\enddocument

\enddocument

\enddocument